\documentclass[reqno,11 pt]{amsart}
\usepackage{amsmath,amssymb,latexsym}
\usepackage[colorlinks=false, linkcolor={black}]{hyperref}
\usepackage{amsthm}
\usepackage{tcolorbox}
\usepackage{color}
\usepackage[all]{xy}
\usepackage{multirow}
\usepackage{physics}
\usepackage{longtable}
\usepackage{caption}
\newcommand{\F}{{\mathbb{F}}}
\newcommand{\Z}{{\mathbb{Z}}}
\newcommand{\Q}{{\mathbb{Q}}}

\newcommand{\N}{{\mathbb{N}}}
\newcommand{\R}{{\mathbb{R}}}

\newcommand{\Gal}{\mathrm{Gal}}

\makeatletter
\newcommand*{\rom}[1]{\expandafter\@slowromancap\romannumeral #1@}
\makeatother

\DeclareFontFamily{U}{wncy}{}
    \DeclareFontShape{U}{wncy}{m}{n}{<->wncyr10}{}
    \DeclareSymbolFont{mcy}{U}{wncy}{m}{n}
    \DeclareMathSymbol{\Sh}{\mathord}{mcy}{"58}

\theoremstyle{plain}

\newtheorem{theorem}{Theorem}[section]
\newtheorem*{theorem'''}{Theorem C}
\newtheorem*{theorem'}{Theorem B}
\newtheorem*{theorem''}{Theorem A}
\newtheorem*{theorem''''}{Proposition D}
\newtheorem{proposition}[theorem]{Proposition}
\newtheorem{rem}[theorem]{Remark}
\newtheorem{lemma}[theorem]{Lemma}
\newtheorem{corollary}[theorem]{Corollary}

\newtheorem{example}[theorem]{Example}

\begin{document}

\title{Binary cubic forms and rational cube sum problem}
%Infinitely many primes in each of the residue classes $1$ and $8$ modulo $9$ are sums of two rational cubes}

\author{Somnath Jha, Dipramit Majumdar \& B. Sury}

\begin{abstract}
{%The classical Diophantine problem of determining  which  integers  can be written as a sum of two rational cubes has a long history; from the earlier works of Sylvester,  Selmer, Satg{\'e}, {Lieman} etc. and up to the recent work of Alp{\"o}ge-Bhargava-Shnidman. 
In this note, we use binary cubic forms to study the rational cube sum problem.  We prove that every non-zero residue class $a \pmod {q}$, for any prime $q$,  contains infinitely many primes which are    sums of two rational cubes. We  also prove (unconditionally) that  for any positive integer $d$, infinitely many primes in each of the residue classes \begin{small}$ 1  \pmod {9d}$\end{small} as well as  \begin{small}$  -1 \pmod {9d}$\end{small},  are sums of two rational cubes.  Further,  for an arbitrary integer $N$, we show there are infinitely many primes $p$ in each of the residue classes \begin{small}$  8 \pmod 9$\end{small} and \begin{small}$1 \pmod 9$\end{small},
such that $Np$ is a sum of two rational cubes.}
\end{abstract}
\maketitle
\section*{Introduction}\label{intro}
%\subsection*{Background and our results:} 
Let us call an integer $n$ a rational cube sum if 
there exist two
rational numbers $a$ and $b$ such that $n=a^3+b^3$. The study of which integers $n$ are rational cube sums  has a rich history and can be traced back to the classical works of Sylvester \cite{syl}, Selmer \cite{sel}, Satg\'e \cite{sa}, {Lieman \cite{lie}} and up to a very recent work of Alp{\"o}ge-Bhargava-Shnidman \cite{abs}.  Without any loss of generality, we may assume that $n$ is cube free and greater than $2$. Then the elliptic curve $Y^3+X^3=nZ^3$, expressed in Weierstrass form  as $E_n:Y^2=X^3-432n^2$, has $E_n(\Q)_{\text{tor}}=0$ and it is easy to see  that $n$ is a rational cube sum $\Leftrightarrow \text{rank}_\Z \ E_n(\Q) > 0$. The important work  of \cite{abs} shows that, when ordered by their absolute value,  a positive proportion of integers are rational cube sums and a positive proportion are not. 

It is a natural to ask which integers primes are  rational cub sums and there has been considerable interest in it, starting from the work in  \cite{syl}. Cube sum problem in congruence classes modulo a prime has also been studied. In an interesting work, Lieman \cite{lie}  using  the analytic properties of $L$-functions of CM elliptic curves together with Coates-Wiles theorem, showed the following:
\begin{theorem}\cite[Theorem 0.1]{lie}\label{lie}
Fix a prime $q>  3$, and a congruence class $a \pmod q$. There exist infinitely many cube-free $n$ congruent to $a \pmod q$ such that $n$ can not be expressed as a sum of two rational cubes.
\end{theorem}

In this article, we  establish a complimentary case to the result of \cite[Theorem 0.1]{lie}:   %  

\begin{theorem''}[Theorem \ref{primesinAP}]
    Let $q$ be a  prime. Each residue class $a \pmod{q}$, for $0< a < q $, contains infinitely many primes which are rational cube sums.
    These primes are  obtained as  values of the cubic form \begin{small}$f(X,Y)=(X+Y)^3-9XY^2$\end{small} at certain  integer points.\qed
\end{theorem''}
The expressibility of  primes as rational cube sums  in the congruence classes modulo $9$ has also been studied in great detail in the literature. There is a conjecture, typically
attributed to Sylvester (see \cite{dv}), which predicts that primes $p \equiv
2,5 \pmod 9$ are not rational cube sums, whereas primes  $p \equiv
4,7,8 \pmod 9$ are rational cube sums. In contrast, primes $p \equiv
1 \pmod 9$ may or may not be rational cube sums. The proof of the
fact that primes $p \equiv 2,5 \pmod 9$ are not cube sums goes back
to the works of P\'epin, Lucas and Sylvester \cite{syl}.
Dasgupta-Voight \cite{dv}, showed that for  primes $p \equiv 4,7
\pmod 9$, both $p$ and $p^2$ are rational cube sums provided $3$ is
not a cube modulo $p$.

The questions on the expressibility of primes  $p \equiv 8 \pmod 9$ are reported to  be  `decidedly more difficult' (see  \cite{dv}, \cite{yi}). It seems the only general known result so far  is due to \cite{yi}, which shows for a prime $p \equiv 8 \pmod 9$, if $x^9 -24x^6 +3x^3 +1 -9(\sqrt[3]{3}-1)x^2(x^3 +1)^2 = 0$ has no solutions in $\F_p$, then at least one of $p$ and $p^2$ is a rational cube sum. 
The works of \cite{dv} and \cite{yi}, although different, use the
theory of (mock) Heegner points. To determine which primes $p \equiv
1 \pmod 9$ are rational cube sums  seems to be a rather subtle
question and was investigated in \cite{rz}. The Birch and
Swinnerton-Dyer (BSD) conjecture predicts that the special value of the
complex $L$-function of $E_p$ at $1$ i.e. $ L(E_p/\Q, 1) = * C_p,$
where $* \neq 0$ and $C_p = 0$ if and only if $p$ is a rational cube
sum. In \cite{rz}, three  efficient methods are given to numerically
test whether $C_p=0$ for a prime $p \equiv 1 \pmod 9$.

In this article, we  show (unconditionally) the existence of infinitely many primes, in each of the congruence   classes $8 \pmod 9$ and $1 \pmod 9$, which are rational cube sums:

\begin{theorem'}[Theorem \ref{primes} \& Lemma \ref{hblem9d}]
For any positive integer $d$ and any integer $a \equiv
\pm{1} \pmod 9$ satisfying $(a,d)=1$, there are infinitely many
primes $p \equiv a \pmod {9d}$ which are rational cube sums.

In particular, infinitely many primes in each of the residue classes $1 \pmod 9$ and $8 \pmod 9$ are rational cube sums. Each of these primes   is a value of the cubic form \begin{small}$f(X,Y)=(X+Y)^3-9XY^2$\end{small} at certain integer point, which depends on the congruence class. \qed
\end{theorem'}

\begin{rem}
In \cite{rz}, the authors gave a list of $22$ primes $p \equiv 1
\pmod{9}$, up to $2000$, for which $C_p =0$ i.e. (assuming BSD) $p$
is a cube sum. Our computation yields this interesting observation
that each of these primes is expressible as $f(x,y)=(x+y)^3-9xy^2$ for some $(x,y)
\in \Q^2$. %For example, $883 = f(\frac{14}{17},\frac{195}{17})$.
Similarly, we have verified that all primes $p \equiv 8 \pmod{9}$,
up to $500$, can be expressed as $f(x,y)$ for some $(x,y) \in \Q^2$.

\end{rem}

 A variant of Theorem B  appears in  Corollary \ref{maincor}. 
 For example, Corollary \ref{maincor} applies to $ 27 \pmod d$ whenever  $3 \nmid d$ but it is  not covered by Theorem B. Further, in Proposition \ref{rank2unconditional}, we exhibit a certain family of primes $p$ which are  congruent to $1 \pmod 9$ and are values of the cubic form  $f(x,y)=(x+y)^3-9xy^2$, such that 
  rank$_\Z \ E_p(\Q)=2$.

For a composite integer $n$ with more than one prime divisor, the results regarding
expressibility of $n$ as a rational cube sum are far more scarce.
After the classical work of Sylvester, Satg\'e \cite{sa}, %\cite{sy}
 others (cf. \cite{jms} and \cite{ms2}) have typically considered integers  of the form $p^iq^j$
where $i,j \leq 2$ and $p,q $ are distinct primes with $p \in \{2,
3, 5\}$  and  $q \equiv 2, 5 \pmod{9}.$
% Using $3$-descent on the Selmer group, some of these results were extended recently in .
Results regarding integers having $3$ or more prime factors seem even scarcer and we are only aware of  results in \cite{cst}
which 
showed that for any odd integer $k \ge 1$, there exist infinitely
many cube-free odd integers $n$ with exactly $k$ distinct prime
factors so that $2n$ is a cube sum. %; similarly, there exist infinitely many cube-free odd integers $n$ with exactly $k$ distinct prime factors such that $2n$ is not a cube sum.
In this context, we have a very general result for an arbitrary integer $N$: 
\begin{theorem'''}[Theorem \ref{Np1}]

       For any integer $N$, there are infinitely many primes $p$, in each of the residue classes $  8 \pmod 9$ and $1 \pmod 9$, so that $Np$ is a rational cube sum.
       Each such prime $p$ is a value of the  cubic form $X^3+NY^3$ at a certain integer point. \qed
 \end{theorem'''}
Moreover,  for $N=\ell$ a  prime (includes $\ell=2,3$), we strengthen
Theorem \ref{Np1} in  Corollary \ref{2pand3q}, to include more
congruence classes, in addition to $\pm 1 \pmod 9$.  Some other infinite families of composite cube sum integers are generated in Corollary \ref{corap}.

Given  any integer $n$, a natural question is whether $n =x^3+y^3$ is solvable   over certain quadratic fields. In Proposition \ref{prop4.2},
we  show that  any $n \in \N$ is a sum of two cubes over the
infinite family of imaginary  quadratic fields $\{\Q
(\sqrt{-3(4nt^3-27)}): t\in \N, t\geq 3\}$.  Using a result of
\cite{mm}, we  discuss a variant of  this result in Prop
\ref{mmresult}. 
%We also apply  Proposition \ref{prop4.2} to generate another infinite family of integers that are rational cube sums.

We show in Proposition \ref{sum of rational cube}  that every integer $n$ can be written (in infinitely many different ways) as a sum of two integers both of which are rational cube sum. For a fixed integer $k \ge 0$, set $\pi_k(X): =\#\{  n \in \mathbb{N}: n \leq  X \text{ and both } n,  n+k \text{ are rational cube sums}\}$. In \cite[Theorem 1.1]{abs}, authors prove that $\pi_0(X) \ge \frac{2}{21}X$ for $X \gg 0$. We have established  (Proposition \ref{densityofpairs}) that
 for any integer $k \ge 1$ and for $X \gg 0$, $\pi_k(X)$ is at least $ \frac{1}{\sqrt[3]{4} \sqrt{k}} X^{\frac{1}{6}}$. % for a constant $c(k)= \frac{1}{\sqrt[3]{4} \sqrt{k}}$.   
%It will be interesting to study the asymptotic behaviour of $\pi_k(X)$ for $k \ge 1$.%that given any integer $k$, there
%are infinitely many integers $n$ such that both $n$ and $n+k$ are rational cube sums.
Finally, we consider the set $C_\Q$ of rational numbers which are sums of two rational cubes. For any non-empty set of places $S$ of $\Q$, we prove in Proposition \ref{strongapprox} that  the image of  $C_\Q$ in $\mathbb{A}_{\Q}^{S}$ is dense, i.e. strong approximation holds for $C_\Q$.

A key ingredient in  all of our results is \cite[Theorem 1]{hbm2}. A celebrated work of  Heath-Brown \cite{hb},
using Sieve theoretic methods, showed that the integer values of the
binary cubic form $X^3+2Y^3$ represents infinitely many primes. This
was generalised by Heath-Brown and Moroz, to a general irreducible
integral binary cubic form in \cite{hbm}, and then  in \cite{hbm2},
in a way which is more amenable to control congruence classes
represented by the primes. In  a previous article \cite{ms} by two of us, using the explicit parametrization of integral points on the curve $X^2+27Y^2=4Z^3$ together with the result in \cite{hbm}, it was shown that there are infinitely many primes $p$,   congruent to either $ 1 \pmod 9$ or $ 8 \pmod 9$, such that $p$ is a sum of two rational cubes. Subsequently, Prof. Moroz and also Prof. Heath-Brown wrote to us pointing out that using  \cite{hbm2}, we can significantly improve our previous result. 

Indeed, in this note, we consider some special integral binary cubic forms (like $(X+Y)^3-9XY^2$,  $XY(X-Y), X^3+NY^3$). Then, we exploit  various polynomial identities satisfied by these special polynomials to express each of their integer values as a sum of two rational cubes. Finally, using infinitude of primes represented by these binary cubic forms in suitable congruence classes \cite{hbm2}, we obtain various infinite families of  rational cube sum integers. Evidently, this approach is different from the recent works based on Selmer groups  of elliptic curves and Heegner points.

%\begin{rem}
%According to \cite{dv},\cite{yi}, \cite{rz} and others, the cube sum
%problem for primes congruent to $8 \pmod 9$  and $1 \pmod 9$, is
%considered to be a difficult and subtle problem. Theorem
%\ref{primes}   is the first unconditional result showing existence
%of infinitely many primes in each of these   classes which are
%rational cube sums.
%\end{rem}

It follows from the results in \cite{hbm2} that for $X\gg 0$, the number of cube sum primes  $p \le X$, obtained in the given congruence class in each of the Theorems A, B and C,  %such that $p$ is rational cube sum and $p \equiv 1 \pmod{9}$ (resp. $p \equiv 8 \pmod{9}$) 
is at least $CX^{\frac{2}{3} - \epsilon}$ for some constant $C>0$ and any $\epsilon >0$ (see Remark \ref{density1}).%, we discuss density of the primes appearing in Theorem \ref{primes}.

\medskip

\noindent {\bf Acknowledgment:} We are very grateful to  Prof.
Heath-Brown and Prof. Moroz and for their insightful comments, suggestions, encouragement and for
answering many questions.% S. Jha is supported by SERB grant CRG/2022/005923.
\section{results}
\subsection{Prime numbers as cube sums: proofs of Theorem A \& Theorem B}%\label{rankcomp}
For the rest of the article, we fix some binary cubic forms: %In this section, we discuss the rational cube sums problem for prime numbers in congruence classes and use \cite[Theorem 1.1]{hbm2} to give the proofs of Theorems A and B. %We begin by considering the following binary forms
\begin{small}
\begin{equation}\label{basicpoly}
    f(X,Y)= (X+Y)^3-9XY^2, \quad g(X,Y) =XY(X-Y) \text{\quad and \quad }f_1(X,Y)= X^3+Y^3-3X^2Y% h(X,Y)= X^2-XY+Y^2.
\end{equation}
\end{small}
%functions $f(X,Y), \ g(X,Y)$ and $h(X,Y)$ are closely related via the following identities.
\begin{proposition}\label{mainresult}
   Put $h(X,Y)= X^2-XY+Y^2$. %and $(e(X,Y), e_1(X,Y))$
    %be the binary forms as defined above.
    %in  \eqref{input1} and  %and \eqref{input2} and
    %recall $h(X,Y)=X^2-XY+Y^2$.  %For $t\in \{f,e\}$,
    We  have the following identity in $\Z[X,Y,Z,W]$:
   \begin{small}
  \begin{equation}\label{id4}
   \begin{split}
    \Big[\big((Zf(X,Y) + 3W f_1(X,Y)\big) +
     \big((Zf_1(X,Y) - W f(X,Y)\big)  \Big]^3  & + \\
    % \text{\qquad \qquad \qquad \qquad}
    \Big[\big((Zf(X,Y) + 3W f_1(X,Y)\big) - \big((Zf_1(X,Y) - Wf(X,Y)\big)  \Big]^3 &
    %\text{\qquad \qquad \qquad \qquad \qquad\qquad \qquad \qquad}
    \\
      = (Z^2+3W^2)\big((Zf(X,Y) + 3W f_1(X,Y)\big)(2h(X,Y))^3. & %\text{\qquad}
   \end{split}
\end{equation}
\end{small}
\end{proposition}

\begin{proof}
    Recall that (cf. \cite{ev}) for a binary cubic form $F(X,Y) = aX^3+ bX^2Y +c XY^2+ dY^3$, one can associate the following invariants: $Disc(F):= D = b^2c^2-4ac^3-4b^3d-27a^2d^2+18abcd$, the quadratic Hessian $H(X,Y) = (b^2-3ac)X^2 + (bc-9ad)XY +(c^2-3bd)Y^2 $ and the cubic covariant $G(X,Y)=(2b^3 +27a^2d-9abc)X^3 +3(b^2c+9abd-6ac^2)X^2Y -3(bc^2 + 9acd - 6b^2d)XY^2 -(2c^3+27ad^2 -9bcd)Y^3$ and these invariants are related by the syzygy
    \begin{small}\begin{equation*}4H(X,Y)^3 = G(X,Y)^2+ 27D F(X,Y)^2. \end{equation*}\end{small}
    The discriminant of $f(X,Y)=3^6$, the quadratic Hessian of $f(X,Y)$ is $3^3h(X,Y)$ and the cubic covariant is $3^5f_1(X,Y)$ and hence the syzygy can be rewritten as $f(X,Y)^2 + 3 f_1(X,Y)^2 = 4 h(X,Y)^3$. From the identities $(a^2+3b^2)(c^2+3d^2)= (ac
+ 3bd)^2 + 3 (ad - bc)^2$ and $(a+b)^3+(a-b)^3 = 2a(a^2+3b^2)$, the identity in
\eqref{id4} can be deduced.
\end{proof}

\begin{corollary}\label{1.2}
  \begin{enumerate}
 
  \item  Putting $Z=1, W=0$ in \eqref{id4}, we obtain
  \begin{small}
        \begin{equation}\label{id1}
            \big(f(X,Y) - 3 g(X,Y) \big)^3 + 27 g(X,Y)^3 = f(X,Y)(X^2-XY+Y^2)^3.
        \end{equation}
    \end{small}
    %\begin{enumerate}
    %    \item 
    %    \item Putting $Z=3, W=-1$, we obtain
    %    \begin{equation}\label{id2}
    %        f(X,Y)^3 - f(Y,X)^3 = 27 g(X,Y)h(X,Y)^3.
    %    \end{equation}
    %\end{enumerate}
    In particular, %t is immediate that
 $f(x,y)$ are rational cube sum  for all   $(x,y) \in \Q^2$.
       \item  Let $n \in \Z$ be expressible as $n = z^2+ 3w^2$ for some $z,w \in \Q$. If there exist $x,y \in \Q$ such that $zf(x,y) + 3w f_1(x,y) =1$, then it follows from \eqref{id4} that $n$ is a cube sum.
     \end{enumerate}
\end{corollary}

\begin{rem}
   The identity \eqref{id4} can be directly verified in {\it SAGE} and alternatively, one can also derive \eqref{id4} using the explicit parametrization of integral points on the curve $X^2 + 27Y^2 = 4Z^3$ (cf. \cite{ms}). Moreover,  as $Disc(f(X,Y)) = 3^6$ is  a sixth power, from the syzygy, one can work out if $f(X,Y)=n$ has a rational solution then $E_{n}: y^2=x^3-432n^2$ has a non-trivial  rational point (cf. \cite[\S 3]{bes}).% Further, since $Disc(u f(X,Y) + 3v f_1(X,Y) ) = 3^6 (u^2+3v^2)^2$, again from the syzygy, one can deduce that if $u f(X,Y) + 3v f_1(X,Y)=1$ has a rational solution then $E_{u^2+3v^2}$ has a non-trivial rational point. In-fact, the conclusion of Corollary \ref{1.2}(1) (resp. Corollary \ref{1.2}(2)) will be valid if we replace $f(X,Y)$ (resp. $uf(X,Y)+3vf_1(X,Y)$) by a binary cubic $F(X,Y)$ whose discriminant is $d^6$ (resp. $d^6n^2$). This point of view is explained in detail in \cite[Theorem 27-29]{bes} using  isogeny induced Selmer groups.  
\end{rem}

Let $T(X,Y) \in \Z[X,Y]$ be an irreducible polynomial  representing
a binary cubic form. Let $a,b,d \in \Z$ such that $\mathfrak
T(X,Y):= T(a+dX,b+dY)$ is a primitive polynomial and assume no prime
divides all the values $\{\mathfrak T(x,y)\mid x, y \in \Z\}$. Then
it follows from \cite[Theorem 1] {hbm2} that $\mathfrak T(x,y)$ attains infinitely
many prime values. Using this, we  prove:

\begin{theorem} \label{primes}
 There are infinitely many primes $p$ in each of the residue classes $  1 \pmod{9}$ and $8 \pmod 9$, such that $p$ is a rational cube sum. Each of these primes $p$ of the form $  1 \pmod{9}$ (respectively $  8 \pmod{9}$) is  a value of the cubic form $f(X,Y)=(X+Y)^3-9XY^2$ at $(-1+3x,-1+3y)$ (respectively  $(1+3x,1+3y)$), for some $x,y \in \Z$.

\end{theorem}

\begin{proof}
Observe  that
$f(X,Y)=(X+Y)^3-9XY^2$ is an irreducible polynomial in $\Z[X,Y]$. Put $\mathfrak F_\pm(X, Y):=f(\mp 1+3X, \mp 1+3Y)$. As $\mathfrak F_\pm(0, 0)=\pm 1$, we see that $\mathfrak F_\pm(X, Y) $ are primitive polynomials in $\Z[X,Y]$ and   no prime divides all the values $\{\mathfrak F_+(x,y)| x, y \in \Z\}$ (respectively $\{\mathfrak F_-(x,y)| x, y \in \Z\}$). Thus by %As $f(1,1)=-1$ and $f(-1,-1)=1$,
 \cite[Theorem 1]{hbm2},  there are infinitely
many primes of the form $\mathfrak F_+(x,y)$ (respectively $\mathfrak F_-(x,y)$)
for $(x,y) \in \Z^2$. Clearly, $\mathfrak F_\pm(x, y) \equiv \pm 1
\pmod{9}$ %and  $\mathfrak F_-(x, y) \equiv 8 \pmod{9} \
for all $ x,y \in \Z$. Now the statements follow from Corollary \ref{1.2}(1).
\end{proof}

Recall, for a prime $p\equiv 1 \pmod 9$, the root number of $E_p: y^2=X^3-432p^2$ is $1$ and by a $3$-descent argument, one can check rank$_\Z \ E_p(\Q) \leq 2$. Thus assuming the parity conjecture, the rank of $E_p(\Q)$ is $0$ or $2$. We now show  (without any assumption on the BSD conjecture or the parity conjecture) that for a certain family of primes $p \equiv 1 \pmod 9$, the Mordell-Weil ranks of $E_p(\Q)$ and $E_{p^2}(\Q)$ equals $2$.
\begin{proposition}\label{rank2unconditional}
Let $p \equiv 1 \pmod{9}$ be a prime and put  $p' \in \{ p, p^2 \}$. Let $(u,v) \in \Z^2$ such that $p'=u^2+3v^2$. Assume both the following conditions hold:
\begin{enumerate}
    \item There exists rational numbers $x,y \in \Q$ such that $f(x,y)=p'$.
    \item There exists rational numbers $z,w \in \Q$ such that $uf(z,w) + 3v f_1(z,w) =1$.
\end{enumerate}
    Then $\mathrm{rk}_{\Z} \  E_{p'}(\Q) =2$ and the $3$-part of the Tate-Shafarevich group  $\Sh(E_{p'}/\Q)[3]$ vanishes.
    \end{proposition}

\begin{proof}
   We prove the result for $p$ and the proof for $p^2$ is similar. We know that \begin{small}$E_p(\Q)_{\text{tor}}=0$\end{small} and  using a $3$-descent argument (for example \cite[\S5]{jms}), one can show that the $\F_3$-dimension of the $3$-Selmer group \begin{small}$S_3(E_p/\Q)$\end{small} of \begin{small}$E_p$\end{small} over $\Q$ is at most $2$. We will show that \begin{small}dim$_{\F_3} \ E_{p}(\Q)/3E_{p}(\Q) \geq 2$\end{small}. Hence  \begin{small}$\mathrm{rk}_{\Z} \  E_{p}(\Q) =2$\end{small} and \begin{small}$\Sh(E_{p}/\Q)[3]=0$\end{small} follows from the descent exact sequence \begin{small}
       $0 \rightarrow \ E_{p}(\Q)/3E_{p}(\Q)  \rightarrow S_3(E_p/\Q) \rightarrow \Sh(E_{p}/\Q)[3] \rightarrow 0.$\end{small}

   Recall that if \begin{small}$X^3+Y^3=p$\end{small}, then \begin{small}$\Big( \frac{12p}{X+Y}, 36p \frac{X-Y}{X+Y} \Big) \in E_{p}(\Q)$.\end{small} By condition $(1)$, it follows  from \eqref{id1} that \begin{small}$P:= \big(12h(x,y), 36f_1(x,y) \big) \in  E_{p}(\Q)$.\end{small} Further using hypothesis $(2)$, from the identity  \eqref{id4}, it follows that \begin{small}$Q:=\big(12ph(z,w), 36p[uf_1(z,w)-vf(z,w)] \big) \in  E_{p}(\Q)$.\end{small}
 Also recall that for a non-trivial rational point \begin{small}$(x,y) \in E_{p}(\Q)$\end{small}, the Kummer map \begin{small}$\delta: E_{p}(\Q) \to \Q(\zeta_3)^*/{\Q(\zeta_3)}^{*3}$\end{small} is given by \begin{small}$\delta\big( (x, y) \big) = y - 12p\sqrt{-3}$\end{small} (for example, see \cite[\S 3]{jms}). It follows that \begin{small}$\delta(P) = \zeta_3^2$ and $\delta(Q) = \zeta_3^2 \pi_p^2 \overline{\pi}_p$\end{small}, where \begin{small}$\zeta = ({-1 + \sqrt{-3}})/{2}$\end{small} and \begin{small}$\pi_p = u + \sqrt{-3}v$.\end{small} As a consequence, it follows that \begin{small}$\delta(P), \delta(Q), \delta(P+Q)$\end{small} and \begin{small}$\delta(P-Q)$\end{small} are all non-trivial elements in \begin{small}$\Q(\zeta_3)^*/{\Q(\zeta_3)}^{*3}$\end{small}. Hence $[P], [Q], [P+Q]$ and $[P-Q]$ are non-trivial in \begin{small}$E_{p}(\Q)/3E_{p}(\Q)$\end{small}, which in turn implies that   \begin{small}dim$_{\F_3} \ E_{p}(\Q)/3E_{p}(\Q)=2$.\end{small}
\end{proof}

\begin{rem}
   % One may wonder how often both the conditions in Prop. \ref{rank2unconditional} hold. In \cite{rz}, assuming BSD, it is shown that  there are $22$ primes $p \equiv 1 \pmod{9}$, $p <2000$ for which $\mathrm{rk}_{\Z} \  E_{p}(\Q)=2$. 
   We have verified via SAGE  that both the conditions of Prop. \ref{rank2unconditional} are satisfied for each prime $p \equiv 1 \pmod 9$ such that $p \leq 2000$ and $p$ is a rational  cube sum.

\end{rem}

\iffalse
\begin{rem}\label{density2}
 Numerical computations carried out by us indicate that number of distinct primes $p$, such that $f(a,b)=p^2$ for $(a,b) \in [-X,X] \times [-X,X]$ is of the
 order $\frac{\sqrt{X}}{\log X}$, and moreover roughly half of these primes are $1$ modulo $9$ and the other half of these primes are $-1$ modulo $9$.
However, to prove that a binary cubic form attains infinitely many
prime square values, one has to obtain an analogue of the lower
bound for the density obtained in the proof of \cite[Theorem
2]{hbm2} using Sieve theory methods. But, we understand following
some correspondences with Prof. Heath-Brown that proving such a
result may be a difficult problem in Sieve theory.
\end{rem}
\fi

From Theorem \ref{primes}, it follows that for $q=2,3$, the non-zero
residue classes modulo $q$  contain infinitely many
primes that are rational cube sums. To establish this result for a general prime $q$, we begin with the following observation for $f(X,Y)=(X+Y)^3-9XY^2$:

\begin{lemma}\label{nonsing}
    Let $q >3$ be a prime and let  $0<a<q$ be an integer. Then $f(X,Y)=a$ is a non-singular plane curve over $\F_q$.
\end{lemma}

\begin{proof}
    The curve is singular over $\F_q$ if and only if  $\exists \ (\alpha,\beta) \in \F_q^2$ such that $f(\alpha,\beta)=a$ and the partial derivatives  $f_X(\alpha,\beta) = f_Y(\alpha,\beta)=0$. First of all,  if  $\alpha \in \F_q$ with $f(\alpha, 0)=a$ i.e. $\alpha^3 =a$, then $f_X(\alpha,0)=f_Y(\alpha, 0) = 3\alpha^2 \neq 0$. Thus, we can assume that $\beta \neq 0$.
    Then \begin{small}\begin{equation}\label{equsing}f_X(\alpha, \beta) = 3 \beta^2( ( \alpha/{\beta} +1)^2 -3) \  \text{ and } \ f_Y(\alpha, \beta) = 3 \beta^2( ( {\alpha}/{\beta} -2)^2 -3).\end{equation}\end{small}Thus it reduces to assume that $3 = \gamma^2$ for some $\gamma \in \F_q$. However, using $q >3$ and $\beta \neq 0$, we get from \eqref{equsing} that % $f_X(\alpha, \beta) = 0 \Rightarrow \frac{\alpha}{\beta}= -1 + \gamma \text{ or } -1 -\gamma$, while $f_Y(\alpha, \beta) = 0 \Rightarrow \frac{\alpha}{\beta}= 2 + \gamma \text{ or } 2 -\gamma$. Since $q>3$, $f_X(\alpha, \beta) =f_Y(\alpha, \beta) = 0$ implies either $-1 -\gamma = 2 + \gamma$ or $-1 +\gamma = 2 - \gamma$,which in turn implies that
  $4 \gamma^2 -9 =0$, which is a contradiction.\end{proof}

\begin{corollary}\label{nonemptyrationalpoints}
    Let $q >3$ be a prime and $0<a<q$ an integer. There exists $(\alpha, \beta) \in \F_q^2$ such that $f(\alpha, \beta) = a$.
\end{corollary}

\begin{proof}
    By Lemma \ref{nonsing}, the cubic curve $f(X,Y)=aZ^3$ is non-singular. If $f(X,Y)-aZ^3$ is reducible, then there is a linear factor  %then its components are curves of genus $0$
    and hence there is a rational point $(\alpha, \beta) \in \F_q^2$ with $f(\alpha, \beta) = a$.
    %hence, each component contains $(q+1)$ rational points (including the point of infinity).
    On the other hand, if $f(X,Y)-aZ^3$ is irreducible,
    then it is a curve of genus $1$, and hence it has $q+1-t$ rational points (including the point of infinity), where $|t| < 2 \sqrt{q} <q$
    (since $q \ge 5$). Thus, in either case,   $\exists \ (\alpha, \beta) \in \F_q^2$ such that $f(\alpha, \beta) = a$.
\end{proof}

Now we are ready to prove Theorem A:% {which is complimentary to the result of \cite{lie} (Theorem \ref{lie}).}

\begin{theorem}\label{primesinAP}
    Let $q >3$ be a prime and $0<a<q$ an integer. Each residue class $a \pmod{q}$ contains infinitely many primes which are rational cube sums.
\end{theorem}

\begin{proof}
    By Corollary \ref{nonemptyrationalpoints}, there exists
    $ (\alpha, \beta) \in \F_q^2$ so that $f(\alpha, \beta)$=$(\alpha +\beta)^3-9\alpha \beta ^2=a$. If $(k_1, k_2) \in \Z^2$ are such that $k_1 \equiv \alpha \pmod{q}$
    and $k_2 \equiv \beta \pmod{q}$, then $f(k_1+qx, k_2+qy) \equiv a \pmod{q}$ for any $(x,y) \in \Z^2$. Recall $f(X,Y) \in \Z[X,Y]$ is  irreducible and  as $(a,q)=1$,  $\mathfrak F(X,Y):= f(k_1+qX, k_2+qY)$ is a primitive polynomial. Further $q \nmid \mathfrak F(0,0)$ and if there is a prime $p \neq q$ dividing all the values $\{\mathfrak F(m,n)\mid m,n \in \Z\}$, then   we can find  $c,d \in \Z$ satisfying $k_1+qc \equiv 1 \pmod p $ and $k_2+qd \equiv 0 \pmod p$.  Then $\mathfrak F(c,d) \equiv 1 \pmod p,$ which is a contradiction. Hence we can apply \cite[Theorem 1]{hbm2} to deduce that there are
    infinitely many primes  of the form $\mathfrak F(x, y)$ for $(x,y) \in \Z^2$ and by construction $\mathfrak F(x, y) \equiv a \pmod q$.  Each of these primes are rational cube sums by Corollary \ref{1.2}(1).
\end{proof}

Following  Theorems \ref{primes} and \ref{primesinAP}, one may ask the following more general question: Let $d$ be a positive integer. When does the residue class $a\pmod{d}$ with $(a,d)=1$, contain infinitely many primes (or contains no primes) which are rational cube sums?

 One can imitate the proof above to show that if there exists $(k_1, k_2) \in \Z^2$ such that $f(k_1, k_2) \equiv a \pmod{d}$, then there are infinitely many primes $p\equiv a \pmod{d}$ which are rational cube sums. Observe  that $f(k,k) = -k^3$ and hence, we deduce:
 \begin{corollary}\label{maincor}
Let $d $ be a positive integer and $(a,d)=1$. If $a$ is a cube modulo $d$, then the residue class $a\pmod{d}$ has infinitely many primes which are rational cube sums. 
\end{corollary}
Theorem B is now immediate from Lemma \ref{hblem9d}, which
was suggested to us by Professor Heath-Brown. Recall,  we have set \begin{small}$f(X,Y) = (X+Y)^3-9XY^2$.\end{small}
%) which extends Theorem \ref{primes} follows from Lemma \ref{hblem9d} and Corollary 1.2 (1). 

%\begin{corollary}\label{hbcor9d}
%Let $d$ be any positive integer. For  any $a \in  \Z$ with  $(a,d)=1$ and $a \equiv
%\pm{1} \pmod 9$,  there are infinitely many
%primes $p \equiv a \pmod {9d}$ which are rational cube sums.
%\end{corollary}

%\begin{proof}
%    The proof is immediate from  the following
%    Lemma \ref{hblem9d} and Proposition \ref{poly1cubesum}. 
%\end{proof}
% Lemma \ref{hblem9d} was communicated to us by Professor Heath-Brown.
\begin{lemma}\label{hblem9d}
 Suppose $a
\equiv \pm{1} \pmod 9$. Let $d$ be any positive integer such that
$(a,d)=1$. Then, there exist integers $r,s$ such that $f(r,s)
\equiv a \pmod {9d}$. Further, there are infinitely many primes $p$
satisfying $p \equiv a \pmod {9d}$ which are values of $f$.
\end{lemma}

\begin{proof}

For the first statement, we prove it for prime power moduli and
apply the Chinese remainder theorem. First, if $p \neq 3$ is a
prime dividing $d$, then $f(r,s) \equiv a \pmod p$ has a solution in
integers $r,s$ by Corollary \ref{nonemptyrationalpoints}.

To obtain integers $x,y$ such that $f(x,y) \equiv a \pmod {p^e}$ for
$e>1$, we apply Hensel's lemma. If $f(r,s) \equiv a \pmod p$, we
show that the partial derivatives $f_x(r,s)$ and $f_y(r,s)$ are not
both $0 \pmod p$. Indeed,
$f_x = 3(x+y)^2-9y^2; f_y = 3(x+y)^2-18xy.$
If $f_x(r,s) \equiv f_y(r,s) \equiv 0 \pmod p$ where $f(r,s) \equiv
a \pmod p$, then we have
\begin{small}
  $$r f_x(r,s) + sf_y(r,s) = 3f(r,s) \equiv 0 \pmod p$$
\end{small}
which is a contradiction.
Therefore, Hensel's lemma applies to give a solution $x,y$ for the
congruence $f(x,y) \equiv a \pmod {p^e}$ for any $e \geq 1$. \\
Next, let $p=3$; then it is easy to see (by induction)  that $f(r,0)=r^3 \equiv a \pmod {3^e}$ has solutions
in integers $r$ whenever $a \equiv \pm{1} \pmod 9$. Observe that, writing  $a =\pm 1+9k, k \in \Z$ it follows  that $(\pm 1+3k)^3 \equiv \pm 1 +9k \equiv a
 \pmod {3^3}$. Now for a general $e$,  assume there exists a solution $r \in \Z$ of $r^3\equiv a \pmod {3^{e-1}}$. Then  using  the $3$-adic expansion of $r$, it is easy   to get a solution $r_1\equiv r \pmod {3^{e-1}}$ such that $r_1^3\equiv a \pmod {3^e}$.  Therefore, the
first assertion is proved.% (we have used corollary 2.7).

Now, if $f(r,s) \equiv a \pmod {9d}$, then consider $\mathfrak F(x,y): = f(r+9dx,s+9dy).$
  $\mathfrak F$ is a primitive polynomial over $\Z$. Also  there is no common prime divisor $p$ of all the
values of $\mathfrak F$. If such a prime exists, then $(p,9d)=1$ because $\mathfrak F(0,0)
\equiv a \pmod {9d}$ and $(a,9d)=1$. But, as before, we can solve for
$x_0,y_0$ such that $r+9dx_0 \equiv 1 \pmod p$ and $s+9dy_0 \equiv
0 \pmod p$. Then, $\mathfrak F(x_0,y_0) \equiv f(1,0) = 1 \pmod p$.\\
%As $\mathfrak F$ is a primitive polynomial over $\Z$ and the values of $\mathfrak F$ has no common prime divisor, 
Hence we  deduce from \cite[Theorem 1.1]{hbm2} that
\begin{small}$\mathfrak F(x,y)$\end{small} takes infinitely many prime values.
\end{proof}

\begin{rem}[Density of primes]\label{density1} %In fact, their result is much stronger, they also provided an asymptotic formula of density of such primes. To elaborate, It is natural
It is evident  that the density of cube sum primes appearing in each of the
cases $\pm1 \pmod 9 $ in Theorem \ref{primes}, is the same as the
density obtained in \cite[Theorem 1]{hbm2}. Let \begin{small}$F \in \{-\mathfrak
F_\pm\}.$\end{small} For \begin{small}$X \gg 0$\end{small},  \cite{hbm2} considers a square \begin{small}$I(X) = (X,
X(1 + \eta)] \times (X, X(1 + \eta)] $\end{small}, where \begin{small}$\eta = (\log X)^{-c}
<1$\end{small} and  \begin{small}$\pi(F,X)$\end{small} counts  the number of primes of the form
\begin{small}$F(a,b)$\end{small} for \begin{small}$(a,b) \in I(X)$\end{small}. By \cite[Theorem 2]{hbm2}, $\pi(F,X)=
\sigma_1(F) \frac{\eta^2 X^2}{3 \log X} \{ 1 + O((\log \log
X)^\frac{-1}{6}) \}$ as \begin{small}$X \to \infty$\end{small}, where \begin{small}$\sigma_1(F)$\end{small} is a positive constant,  depending on $F$. 

It follows that for $X\gg 0$, number of cube sum primes $p \le X$, obtained  in the given congruence class in each of Theorems A, B and C,  is at least $cX^{\frac{2}{3} - \epsilon}$ for any $\epsilon >0$ and for some  constant $c >0$.
\end{rem}

%%%%%%%%%%%%%%%%%%%%%%%%%%%%%%%%%%%%%%%%%%%%%%%%%%%%%%%%%%% New Section
%\section{Further polynomial identities and cube sum composite numbers }\label{sec3}
%In this section, we generate several infinite  families of composite
%integers that are rational cube sums. 

\subsection{Composite numbers as cube sums: proof of Theorem C}
We will use  the following reformulation of the cube sum problem.  Although it seems to be known to experts, we write it down for the sake of completeness. %using the identity
%\eqref{id2}.
%Note that yet another reformulation appears in Prop. \ref{obvious}. 

\begin{lemma} \label{firstequivalence}
   Any  integer $n $ of the form $n=xy(x-y)$ with $x, y \in \Q$ is a rational cube sum.  Conversely, if a cube free integer $n \geq 2 $ 
 is a rational cube sum then there exists $ x,y \in \Q$  such that  $n=xy(x-y)$.
\end{lemma}
\begin{proof}

   The first assertion is immediate from    \eqref{id4} by setting  \begin{small}$Z=3, W=-1$. \end{small}
The converse can be deduced using  the fact that \begin{small}$E_{n}$\end{small} is
$3$-isogenous over \begin{small}$\Q$\end{small} to the curve \begin{small}$ Y^2=X^3+16n^2$.\end{small}  \end{proof}
\iffalse
\begin{rem}\label{cubesumzeta3}
 In fact, as   
$ \mathrm{rk}_\Z \ E_{n}(\Q(\zeta_3)) = 2 \ \mathrm{rk}_\Z \ E_{n}(\Q)$,
we can also deduce:\\ $n\in \Z$ is a rational cube sum
$\Leftrightarrow \exists \ x, y \in \Q(\zeta_3)$ such that
$n=xy(x-y)$.
\end{rem}
\fi
%\begin{example}
%It is perhaps not obvious to express even a small number like $42$
%as a rational cube sum; using our identity \eqref{id2}, we deduce
%\begin{small}
%$g(-6,1)= 42= ( {449}/{129} )^3 + ( {-71}/{129})^3$.\end{small}
%\end{example}

Using Lemma \ref{firstequivalence} and \cite[Theorem 1]{hbm2}, we are now ready to prove Theorem C.
\begin{theorem}\label{Np1}
    For any integer $N$ and for any $(x,y)\in \Z^2$ with $xy \neq 0$, the integer $N(x^3+Ny^3)$ is a rational cube sum. In particular, for every $N\in \N$, there are infinitely many primes $p$ in each of the residue classes $1\pmod 9$ and  $8 \pmod 9$ such that $Np$ is a rational cube sum.
    \end{theorem}

   \begin{proof} Consider the binary cubic form   $X^3+NY^3 $ and observe the identity $N(x^3+Ny^3)=\frac{1}{x^3}g(-Ny^3,x^3)$. We may assume that $N$ is not a cube   (by Theorem \ref{primes}). Then applying  \cite[Theorem 1]{hbm2}, the result follows (by a similar argument as in Theorem \ref{primes}).
   \end{proof}

   The cube sum problem for $2p$ and  $3p$, for a prime $p$ are also  stated as cases of Sylvester's Conjecture in the literature (cf. \cite{syl}, \cite{sa}, %\cite{sy}, 
   \cite[Corollary 5.9]{jms} for some relevant works). %Note that in Theorem \ref{Np1} and Corollary \ref{2pand3q}, we had $p=a^3+Nb^3$, for some $a,b \in \Z$ i.e. $N$ is a cube in $\mathbb F_p$. On the other hand, for any prime $p\equiv 1,7 \pmod 9$ so that $2$ is not a cube in $\mathbb F_p$, it is shown in \cite[Corollary 5.9]{jms} that $2p$ is not a rational cube sum. For any prime $p\equiv 2 \pmod 9$, it was shown in \cite{sa} that $2p$ is a cube sum. However, if $p\equiv 5 \pmod 9$, then $2p$ is a not cube sum \cite{syl}. Recently, it was shown in \cite{sy} that for any prime $p\equiv 2,5 \pmod 9$,  $3p$ and also $3p^2$ are   cube sums. 
   We strengthen
   Theorem \ref{Np1} in the special case when $N$ is a prime: %for $N=2,3$ and, in fact,  more generally for $N=\ell$ for a prime $\ell$. One can prove similar results for $N=\ell^2$. %The proofs are   omitted.
\begin{corollary}\label{2pand3q}
    \begin{enumerate}
    %\item For each $a\in \{1,7,8\}$, there are infinitely many primes $p \equiv a \pmod 9$ such that $2p$ is a rational cube sum.
    \item For each $b\in \{1,2,7,8\}$, there are infinitely many primes $p \equiv b \pmod 9$ such that $2p$ is a rational cube sum.
     \item For each $b\in \{1,4,7,8\}$, there are infinitely many primes $p \equiv b \pmod 9$ such that $3p$ is a rational cube sum.

     \item  Let $\ell$  be a fixed prime and so that $\ell  \equiv a \pmod{9}$ with $a \in \{ 1,2,7,8 \}$ (respectively $a \in \{4, 5\}$).  Then there are infinitely many primes $p \equiv b \pmod 9$ where $b\in \{1,2,7,8\}$ (respectively $b\in \{1,4,5,8\}$) such that $\ell p$ is a rational cube sum.
     
     %\item Let  $a \in \{ 1,4,5,8 \}$. Let $\ell$ be a fixed prime such that $ \ell \equiv a \mod{9}$. For each $b\in \{1,2,7,8\}$, there are infinitely many primes $p \equiv b \pmod 9$ such that $\ell^2 p$ is a rational cube sum.
     %\item Let $\ell \equiv a \mod{9}$ be a fixed prime, where $a \in \{ 2,7 \}$. For each $b\in \{1,4,5,8\}$, there are infinitely many primes $p \equiv b \pmod 9$ such that $\ell^2 p$ is a rational cube sum.
    \end{enumerate}

\end{corollary}

Next, we draw several other corollaries of the Lemma \ref{firstequivalence}. %We state them separately simply because each has a different flavour.
\begin{corollary}\label{corap}
  \begin{enumerate}
  \item For $a, b \in \Z$, we can write $g(a, b)= 2 (a-b) \frac{a}{2} b$. Conversely, given $a, d \in \Q$, we have $2a(a+d)(a+2d)=g(2a+2d, a+2d)$. Thus a cube-free integer $n>1$ is a rational cube sum  $\Leftrightarrow n$ is twice the product of three rational numbers in AP. In particular,  $2n(n+k)(n+2k)$ is a rational cube sum $\forall \ n,  k \in \Z$.

         \item  For $x,y \in \Z, \ xy(x^{3k+1} \pm y^{3k+1})=g(\mp x^{3k+1}, y^{3k+1})$  are rational cube sums. In particular,  the product of any two consecutive integers is a rational cube sum. 
        
         % \item Given any three integers $a,b,c$, the integer $(a-b)(b-c)(c-a)=g(b-a,b-c)$ is a rational cube sum.
        % \item For any $n, k \in \Z$, $n(n+k)(2n+k)=g(-n,n+k)$ is a rational cube sum. In particular,  for any  $n \in \Z$,  $2n(n+2)$ is a rational cube sum.
        % \item For any $d, k \in \Z$,  $2d(d^2-k^2)=g(-d-k,d-k)$ is a rational cube sum.
       % \item  If $(x,y)\in \Z^2$ satisfies the Pell's equation $X^2-NY^2=1$, then $Nx^2y^2=g(x^2,1)$ is a rational cube sum.
       \item %If $ABC$ is a rational right angle triangle  whose sides have lengths $a,b,c \in \Q$ respectively,
       For a Pythagorean triplet $(a,b,c)$, $g(c^2,a^2)=(abc)^2$ is a rational cube sum. So there are infinite many perfect square integers which are cube sums.
       \item %Let $F_n$ denote the $n^{th}$ term of the Fibonacci Sequence. %(resp. Lucas number); the convention here is $F_0=0, F_1=1, L_0=2,L_1=1$.Then a
       Product of the three consecutive terms in the Fibonacci sequence is a  cube sum.

     %  \item As a special case of Corollary \ref{corap}, $2n(n+1)(n+2)=g(2n+2,n+2)$ is a rational cube sum. It is the area of an isosceles  Heron triangle whose sides have lengths $n^2+2n+2, n^2+2n+2$ and $2(n^2+2n)$.

%         \item Integers of the form $n=d(8d+1), d \in \Z$  are rational cube sums.  In particular, for any integer $k$, $ \frac{1}{2}k(k+1)(2k+1)^2$ is a rational cube sum.
%        \item For any integer $n$,  $n^2-16$, $3n^2+16$ and $n^2+432$ are rational cube sums.
  %     More generally, integers of the form $(4T^3)2+3V^2$ and $\frac{1}{4}(L^2+ 27 (M^3)^2)$ are rational cube  sums.
   %     \item If the  polynomial $T(X)=aX^2 -bX+c \in \Z[X]$ has a root in $\Q(\zeta_3)$, then $abc$ is a rational cube sum.
          \end{enumerate}
\end{corollary}

\subsection{Cube sums over imaginary quadratic fields}\label{sec4}

We say $n\in \N$ is a  cube sum over a number field $F$ if $n=x^3+y^3$ with $x,y \in F$.
We discuss expressibility  of  an integer  as a cube sum over imaginary quadratic fields. Firstly, we show, under an assumption on the Tate-Shafarevich  group $\Sh(E_{n}/\Q)$ of $E_n: y^2=x^3-432n^2$, that every $n \in \N $ is a cube sum over infinitely many imaginary quadratic fields.

%Recall, $E_n$ denotes  the elliptic curve   $y^2 = x^3 - 432n^2$.
%For any prime $p$, let $\Sh(E_n/\Q)[p^\infty] $ denote the
%$p$-primary torsion part of the Tate-Shafarevich  group of $E_n/\Q$
%and let $S_{p^\infty}(E_n/\Q)$ be the $p^\infty$-Selmer group of$E_n/\Q$. Firstly, we show, under an assumption on
%$\Sh(E_n/\Q)[p^\infty]$, that $n \in \N $ is a cube sum over
%infinitely many imaginary quadratic fields.
\begin{proposition}\label{mmresult}
Let $n$ be any positive integer. Assume that the $p$-primary part of $\Sh(E_{n}/\Q)$ is finite for some prime $p \geq 5$. Then there are infinitely many  imaginary quadratic  fields $\Q(\sqrt{-D_n})$ such that $n $ is a  cube sum over $\Q(\sqrt{-D_n})$, for each $D_n$.
\end{proposition}
\begin{proof}
 We may assume that $n$ is not a cube sum over $\Q$. (In particular,  $n>2$.)  Then $E_n(\Q)$ is finite.
 (By assumption,) choose a prime $p > 3$ such that  $\Sh(E_n/\Q)[p^\infty]$ is  finite. Then the $p^\infty$-Selmer group $S_{p^\infty}(E_n/\Q)$ is finite.
 Now  $E_n$ has CM by $\Q(\zeta_3)$ and $p \nmid \#O^\times_{\Q(\zeta_3)}$, it follows from Rubin's work \cite{ru}  on the Iwasawa main conjecture for imaginary quadratic fields
 that the complex $L$-value $L(E_n/\Q,1) \neq 0$. Further, by results of  Bump-Friedberg-Hoffstein, or Murty-Murty \cite[Corollary to Theorem 2]{mm},
 there are infinitely many imaginary quadratic fields $\Q(\sqrt{-D_n})$ such that for  the quadratic twists $E_n^{D_n}$ of $E_n$, the $L$-functions
 $L(E_n^{D_n}/\Q,s)$ have a simple zero at $s=1$.  Then it is known, by the Gross-Zagier theorem together with  Kolyvagin's or  Rubin's %\cite{ru2}
result that for  any such $D_n$, the Mordell-Weil rank of
$E^{D_n}_n(\Q)$   is $1$. Thus we get that $\text{rank}_\Z \
E_n(\Q(\sqrt{-D_n}))=1$  as well and consequently $n$ is a  cube sum
over $\Q(\sqrt{-D_n})$ for each of those   $D_n$'s.
\end{proof}
Now we give  an explicit and unconditional construction of infinitely many imaginary quadratic
fields over which $n$ can be expressed as a  cube sum.
\begin{proposition}\label{prop4.2}
Let $n$ be any positive integer. For every integer $t {\geq 3} $,  $n $ is a  cube sum over the   imaginary quadratic  field $K_{n,t}=\mathbb Q \big(\sqrt{-3(4nt^3-27)}\big)$.
\end{proposition}
\begin{proof}
%We construct an explicit binary cubic form and then make use of  the
%strategy outlined by  Mordell, Evertse (cf. \cite[\S 3]{ev}). For
%any binary cubic form $F(X,Y)$, recall that  the `quadratic
%covariant' $H(x,y):= -\frac{1}{4}
%\bigg(\frac{\partial^2{F}}{\partial X^2}
%\frac{\partial^2{F}}{\partial Y^2} - (\frac{\partial^2{F}}{\partial
%X \partial Y})^2 \bigg)$ has the discriminant $=-3D$, where $D$ is
%the discriminant of $F$. Further, if we set the `cubic covariant'
%$G(x,y):= \frac{\partial F}{\partial X} \frac{\partial H}{\partial
%^Y} - \frac{\partial F}{\partial Y} \frac{\partial H}{\partial X}$,
%then one has $4H(X,Y)^3=G(X,Y)^2+27DF(X,Y)^2.$

For any positive integer  $t$, consider the binary cubic form $F_{n,t}(X, Y)=(X+Y)^3-nt^3XY^2$. It is easy to see that  $F_{n,t}(X, Y)  \in \Z[X,Y]$ is an  irreducible polynomial  for $t>2$ with discriminant $D=D(F_{n,t})=(nt^3)^2(4t^3n-27)$ and thus for any given $n$ and  for every choice of $t>2$, the discriminant $D(F_{n,t})$ is positive. Let $H_{n,t}(X,Y)$ and  $G_{n,t}(X,Y)$ respectively  denote the quadratic Hessian and  cubic covariant of $F_{n,t}(X,Y)$ (see Prop. \ref{mainresult}). Setting $K_{n,t}:=\mathbb Q (\sqrt{-3D})=\mathbb Q (\sqrt{-3(4nt^3-27)})$, it follows that $F_{n,t}(X,Y) $ is an irreducible polynomial in $K_{n,t}[X,Y]$ as well \cite[Page 122]{ev}. Further, put $U^\pm_{n,t}(X,Y)=\frac{1}{2}(G_{n,t}(X,Y)\pm 3\sqrt{-3D}F_{n,t}(X,Y)) $. %and $V(x, y)= \frac{1}{2}(G(x,y)-3\sqrt{-3D}F(X,Y)) $,
Then we get $U^+_{n,t}(X,Y)U^-_{n,t}(X,Y)=H_{n,t}(X,Y)^3$ \cite[Eq.
10]{ev}. As $U^+_{n,t}(X,Y)$ and $ U^-_{n,t}(X,Y)$ have no common
factors, it follows that each of them is a cube of some homogeneous
linear forms $\xi^\pm_{n,t}(X,Y) \in K_{n,t}[X,Y]$, respectively.
One obtains that $G_{n,t}(X,Y)= \xi^+_{n,t}(X,Y)^3+\xi^-_{n,t}(X,
Y)^3$ \cite[Eq. 11]{ev}.

Observe that $G_{n,t}(1,0)=-27nt^3$. Thus any given $n \in \N$ is a  cube sum over the imaginary quadratic field $\Q (\sqrt{-3(4nt^3-27)})$ for every choice of $t\geq 3$.
\end{proof}
\begin{rem}
Note that, for any positive integer  $n$, $K_{n,t}=\mathbb Q \big(\sqrt{-3(4nt^3-27)}\big)$, as we vary $t$ in $\N_{\geq 3}:=\{n \in \N \mid n \geq 3\}$,  represents  infinitely many imaginary quadratic fields. First, observe there are infinitely many
 prime divisors of the values of $D(t): = 4nt^3 - 27$ as $t$ varies in $\N_{\geq 3}$ %$S= \{ m \in \N \mid m >2 \} $
(If $p_1=3, p_2, \dots, p_r$ are the only  prime divisors, then we arrive at a contradiction by considering $D(3p_1\cdots p_r)$.) Now, for any prime $p \ge 5$, if $p \mid D(t)$ for some $t \in \N_{\geq 3}$, then $D(t+p) \equiv D(t) + 12nt^2p \pmod{ p^2}$; thus although $p$ divides both $D(t)$ and $D(t+p)$, $p^2$ can not divide both $D(t)$ and $D(t+p)$, which in-turn implies that infinitely many primes occur in the factorisation of the square-free part of $D(t)$ as $t$ varies in $\N_{\geq 3}$.
\end{rem}
\iffalse
\begin{corollary}\label{expl-quadratic}
Let $T(X,Y)\in \Z[X,Y]$ be an irreducible binary cubic such that the discriminant $D $ of $T(X,Y)$ satisfied $D=n^2$ for some $n \in \N$. Then $nT(x,y)$ is a rational cube sum  for every $x,y \in \Q$. In particular, every $m \in  \{k^2+k+7 \mid k\in \Z\}$ is a rational cube sum.
\end{corollary}
\begin{proof}
Following the proof of Prop. \ref{prop4.2}, in this setting, we can express $3\sqrt{-3D}T(X,Y)= \xi^+(X,Y)^3+\xi^-(X, Y)^3$, for some $\xi^\pm_{n,t}(X,Y) \in \Q(\sqrt{-3D})[X,Y]$. Now putting $D=n^2$, it is easy to see that for every $x,y \in \Q$, $nf(x,y)$ is a cube sum in $\Q(\zeta_3)$ and hence over $\Q$.

For the second part, for any  $k \in \Z$, consider the irreducible cubic polynomial $T_k(X,Y)=X^3 - (k-1)X^2Y - (k+2)XY^2 -Y^3.$ Then the discriminant $D(T_k)$ equals $ (k^2+k+7)^2$. Hence for every $x,y \in \Q$, $(k^2+k+7)F_k(x,y)$ is a rational cube sum. In particular, taking $x=1$ and $y=0$, we obtain that integers of the form $k^2+k+7$ are rational cube sums.
\end{proof}

\begin{rem}
    Let $T(X) \in \Z[X]$ be a monic irreducible cubic polynomial such that $\Gal(T) \cong \Z/3\Z$. Then from  Corollary \ref{expl-quadratic}, it follows that $\sqrt{Disc(T)}$ is rational cube sum. An explicit parametrization of all monic irreducible trinomials $T(X)=X^3-aX+b$ with $\Gal(T) \cong \Z/3\Z$ can be found in \cite[Theorem 4.6]{ms}.
\end{rem}
\fi

\subsection{Results on the set of cube sums}
%We first consider the set of rational cube sums. 
Let $C$ (resp. $C_\Q$) be the set of integers (resp. rational numbers) which are rational cube sum.% i.e. $C:= \{ q \in \Q \mid q = x^3 +y^3 \text{ for some } x,y \in \Q \}$ and put $C(\Z):=C \cap \Z$.

\begin{proposition}
\label{sum of rational cube}
Given any integer $k$, there are infinitely many integers $n$ such that both $n$ and $n+k$ are rational cube sums. In-particular, $C + C= \Z$.
\end{proposition}

\begin{proof}
We may assume  $k \neq 0$.  For each $i \geq 1$, $g\big( \frac{1}{2i} + 2i^2k, \frac{1}{2i} - 2i^2k \big) =  k - 16i^6k^3$ is a rational cube sum. %On the other hand, $16t^6k^3 =  (2t^2k)^3 +  (2t^2k)^3$ is clearly a cube sum.
Thus $n_i(k) := -16i^6k^3$ and $k+ n_i(k)$ are rational cube sums.
\end{proof}

 Observe that if  $n^6 \le \frac{X}{16k^3}$, then  $n^6(16k^3) - k <X$ and  both $16n^6k^3 - k$ and $16n^6k^3$ are rational cube sums. Thus we  deduce Proposition \ref{densityofpairs} from Proposition \ref{sum of rational cube}:

 \begin{proposition}\label{densityofpairs}
    For any integer $k \ge 1$, $\pi_k(X)$ is at least $ \frac{1}{\sqrt[3]{4} \sqrt{k}} X^{\frac{1}{6}}$ for $X \gg 0$.% for a constant $c(k)= \frac{1}{\sqrt[3]{4} \sqrt{k}}$.   
\end{proposition}

 We remark that for specific values of $k$, the lower bound of $\pi_k(X)$ can be improved, for example, since both $n^2 -n$ and $n^2-n+7$ are rational cube sums, $\pi_7(X) \ge \sqrt{X}$ for $X \gg 0$. \\

%\begin{rem}
%For a fixed integer $k \ge 0$, let $\pi_k(X): =\#\{  n \in \Z: |n| \leq  X \text{ and both } n,  n+k \text{ are rational cube sums}\}$. Then \cite[Theorem 1.1]{abs} implies that $\pi_0(X) = O(X)$. It will be interesting to study the asymptotic behaviour of $\pi_k(X)$ for $k \ge 1$.
%\end{rem}
%\vspace{2mm}
One can show that $C_\Q$ is dense in $\R$ and $\Q_p$ for all prime $p$. In fact, strong approximation holds for $C_\Q$. The key idea behind the proof is the simple observation that given any two non-zero rational numbers  $r,s$, the rational number $ r- \frac{s^3}{r} = g(\frac{r}{s}, \frac{s^2}{r})$ is a cube sum.

\begin{proposition}\label{strongapprox}
    Let $S $ be a non-empty set of places of $\Q$. The image of $C_\Q$ in  $\mathbb A_\Q^S =\prod_{v \notin S}^{\prime} \Q_v$, via the diagonal embedding $C_\Q \hookrightarrow \mathbb A_\Q^S$, is dense.
\end{proposition}

\begin{proof}
     Given any finite set $T $ of places of $\Q$ disjoint from $S$, elements $x_v \in \Q_v$ where $v \in T$ and  for any $\epsilon >0$, we have to show that: there exists $y \in C_\Q$ such that $|y - x_v|_v < \epsilon$ for all $v \in T$ and $|y|_v  \le 1$ for all $v \notin (S \cup T) \setminus \{ \infty \}$. It suffices to establish the result  when $S$ is singleton. We will write down the proof for the case $S = \{ \ell \}$, a finite place and the proof for the case $S = \{ \infty \}$ is similar.\\
    Further, without any loss of generality, we may assume that  $T = \{ p_1, \dots, p_r, \infty \}$.  Now we are given elements $x_i \in \Q_{p_i}$ for $1 \le i \le r$, $x_{r+1} \in \R$ and an arbitrary $\epsilon$ with $0< \epsilon <1$. By strong approximation for $\Q$, there exists $x \in \Q$ such that $|x - x_i|_{p_i} = p_i^{-n_i} < \frac{\epsilon}{2}$ for  $1 \le i \le r$, $|x - x_{r+1}| < \frac{\epsilon}{2}$ and $x \in \Z_p$ for all $p \notin \{ \ell, p_1, \dots, p_r \}$. We fix such an  $x$ and  assume that $x \notin C_\Q$ (otherwise, we are done).  Write $x = \pm \ell^{m}p_1^{m_1} \cdots p_r^{m_r} q_1^{t_1} \cdots q_s^{t_s}$ with $m, m_1, \dots, m_r \in \Z$, $t_1, \dots, t_s \in \N$ and $q_i$'s are integer primes disjoint from $S\cup T$. We denote  $ max\{t_1, \dots, t_s, n_1, \dots, n_r, |m_1|, \dots, |m_r|, |m| \} $ by $N$ and  set $w:= (p_1 \cdots p_r q_1 \cdots q_s)^{N}$. We fix an integer $M$ such that $\ell^{-3M} < \frac{\epsilon |x|}{2w^3 }$ and define $y: = g( \frac{x\ell^M}{w}, \frac{w^2}{x\ell^{2M}})=x - \frac{w^3}{x\ell^{3M}}$. Then $y \in C_\Q$ and moreover, as $\frac{w^3}{x\ell^{3M}} \in \Z[\frac{1}{\ell}]$, it follows that, $y \in \Z_p$ for all $p \notin \{ \ell, p_1, \dots, p_r \}$. Finally, note that, $|y -x_i|_{p_i} = |x -x_i|_{p_i} < \epsilon$ and $|y -x_{r+1}| \le |x - x_{r+1}| + \frac{w^3}{\ell^{3M}|x|} < \epsilon$.   %A simpler proof works if $\infty \not\in T$.
\end{proof}

\medskip

\begin{small}
\noindent{\bf S. Jha}, Department of Mathematics \& Statistics, IIT Kanpur, Kanpur 208016, India (jhasom@iitk.ac.in)

\medskip

\noindent {\bf D. Majumdar}, Department of Mathematics, IIT Madras,
Chennai 600036, India (dipramit@gmail.com)

\medskip

\noindent {\bf B. Sury}, Stat-Math Unit, Indian Statistical
Institute, 8th Mile Mysore Road, Bangalore 560059, India
(surybang@gmail.com)

\end{small}
\end{document}